\documentclass[11pt]{article}
\usepackage{amsmath,amsfonts,amssymb,amsthm}

\title{Orbifolds are not commutative geometries}

\author{Adam Rennie\dag
\thanks{email: \texttt{rennie@math.ku.dk},
\texttt{varilly@cariari.ucr.ac.cr}}
\word{and}%
Joseph C. V\'arilly\ddag$\Vert^*$ \\[6pt]
\dag Institute for Mathematical Sciences,
University of Copenhagen\\
Universitetsparken 5, DK-2100 Copenhagen, Denmark\\[6pt]
\ddag Departamento de F\'{\i}sica Te\'orica I,
Universidad Complutense,\\
Madrid 28040, Spain\\[6pt]
$\Vert$Departamento de Matem\'aticas,
Universidad de Costa Rica, \\
2060 San Jos\'e, Costa Rica}

\date{21 March 2007}

\advance\topmargin by -\headheight
\advance\topmargin by -\headsep
\evensidemargin=\oddsidemargin
\makeatletter
\def\section{\@startsection{section}{1}{\z@}{-3.5ex plus -1ex minus
  -.2ex}{2.3ex plus .2ex}{\large\bf}}
\def\subsection{\@startsection{subsection}{2}{\z@}{-3.25ex plus -1ex
  minus -.2ex}{1.5ex plus .2ex}{\normalsize\bf}}
\makeatother

\numberwithin{equation}{section} 

\theoremstyle{plain} 

\theoremstyle{definition} 

\DeclareMathOperator{\Aut}{Aut}   
\DeclareMathOperator{\Cliff}{{\C\ell}} 
\DeclareMathOperator{\diag}{diag} 
\DeclareMathOperator{\Dom}{Dom}   
\DeclareMathOperator{\End}{End}   
\DeclareMathOperator{\Id}{Id}     
\DeclareMathOperator{\spec}{sp}   
\DeclareMathOperator{\Tr}{Tr}     
\DeclareMathOperator{\tsum}{{\textstyle\sum}} 
\DeclareMathOperator{\vol}{vol}   

\newcommand{\A}{\mathcal{A}}  
\newcommand{\B}{\mathcal{B}}  
\newcommand{\al}{\alpha}      
\newcommand{\C}{\mathbb{C}}   
\newcommand{\cc}{\mathbf{c}}  
\newcommand{\CDA}{\mathcal{C_D(A)}} 
\newcommand{\Coo}{C^\infty}   
\newcommand{\D}{\mathcal{D}}  
\newcommand{\dl}{\delta}      
\newcommand{\Dslash}{{D\mkern-11.5mu/\,}} 
\newcommand{\Ga}{\Gamma}      
\renewcommand{\H}{\mathcal{H}}  
\newcommand{\half}{\tfrac{1}{2}} 
\renewcommand{\L}{\mathcal{L}} 
\newcommand{\La}{\Lambda}     
\newcommand{\la}{\lambda}     

\newcommand{\ox}{\otimes}     
\newcommand{\oxyox}{\otimes\cdots\otimes} 
\newcommand{\R}{\mathbb{R}}   
\newcommand{\rt}{\triangleleft} 
\newcommand{\sg}{\sigma}      
\newcommand{\Trw}{\Tr_\omega} 
\newcommand{\Z}{\mathbb{Z}}   
\newcommand{\8}{\bullet}      
\renewcommand{\.}{\cdot}      

\def\<#1>{\langle#1\rangle} 

\newcommand{\row}[3]{{#1}_{#2},\dots,{#1}_{#3}} 
\newcommand{\set}[1]{\{\,#1\,\}}  

\newcommand{\word}[1]{\quad\mbox{#1}\quad} 

\newcommand{\comment}[1]{\relax}  

\begin{document}

\maketitle

\begin{abstract}
In this note we show that the crucial orientation condition for
commutative geometries fails for the natural spectral triple of an
orbifold $M/G$.
\end{abstract}

\section{Introduction} 

Recently, in \cite{Crux}, it was shown that a spectral triple
$(\A,\H,\D)$ with $\A$ commutative and unital is the Dirac triple of a
spin$^c$ manifold $(C^\infty(M),L^2(M,S),\Dslash)$ if and only if
$(\A,\H,\D)$ satisfies certain additional conditions. These conditions
are a mild strengthening and generalisation of those proposed
in~\cite{ConnesGrav}. (A spectral triple with $\A$ commutative and
satisfying some version of these conditions could thus be called a
commutative geometry, hence the title).

In this note we show precisely why orbifolds of the form $M/G$, with
$G$ a finite group acting by isometries on $M$, fail to satisfy these
conditions. For such orbifolds we can define a commutative spectral
triple, and though we do not show it here, these spectral triples do
satisfy the majority of those conditions. However, the crucial
orientation condition, which provides a volume form, tangent bundle
and ultimately gives the local coordinates, fails for the spectral
triples of the orbifolds.

After a brief recollection of the conditions, we go directly to the
demonstration of the failure of the orientation condition. To finish,
we briefly comment on more general orbifolds.

\section{Geometric conditions for commutative spectral triples}
\label{sec:conds}

In \cite{Crux} we have shown that under a small number of postulates,
mildly extending those of \cite{ConnesGrav}, any spectral triple
$(\A,\H,\D)$ whose algebra $\A$ is \textit{commutative} and unital
satisfies $\A = \Coo(M)$ where $M$ is a smooth closed manifold
carrying a spin structure, $\H$ is its $L^2$-spinor space and $D$ is a
Dirac operator plus an endomorphism of the spinor bundle.

Recall that a (unital) spectral triple consists of a unital algebra
$\A$ faithfully represented on a Hilbert space~$\H$, and a selfadjoint
operator $\D$ on~$\H$, with compact resolvent, such that $\A$
preserves its domain $\Dom\D$ and $[\D,a]$ extends to a bounded
operator on~$\H$ for each $a \in \A$.

Assume furthermore that $\A$ is commutative. Denote by $A$ its norm
completion in $\B(\H)$ and by $M = \spec(A) = \spec(\A)$ its character
space (Gelfand spectrum), which is assumed to be a separable
$C^*$-algebra. The following postulates, of an algebraic or
operatorial nature, are needed for the reconstruction of the manifold,
spin structure and Riemannian metric.

\begin{enumerate}
\item
\textit{Dimension}:
For some $p \in \{1,2,3,\dots\}$, the operator
$\<\D>^{-1} := (1 + \D^2)^{-1/2}$ lies in $\L^{p,\infty}$ and 
$\<\D>^{-p} \in \L^{1,\infty}$ has positive Dixmier traces:
$\Trw \<\D>^{-p} > 0$ for all~$\omega$.

When $p$ is even, $\H$ is $\Z_2$-graded by a selfadjoint unitary $\Ga$
commuting with $\A$ and anticommuting with~$\D$. When $p$ is odd,
$\H$ is ungraded and we put $\Ga = 1$ for convenience.

\item
\textit{Regularity}:
$(\A,\H,D)$ is a $QC^\infty$ in the sense of \cite{CareyPRS}: if
$\dl(x) := [|\D|,x]$ for $x \in \B(\H)$, then 
$\A \cup [\D,\A] \subseteq \bigcap_{m=1}^\infty \Dom \dl^m$.
Moreover, $\A$ is complete in the locally convex topology given by the
seminorms $q_m(a) := \|\dl^m a\|$ and $q'_m(a) := \|\dl^m([\D,a])\|$,
for $m = 0,1,2,\dots$.

\item
\textit{Finiteness}:
The prehilbert space $\H_\infty = \bigcap_{m\geq 1} \Dom\D^m$ is a
finitely generated projective $\A$-module.

\item
\textit{Absolute continuity}: If $a > 0$ in $\A$, then
$\Trw a\<\D>^{-p} > 0$ for all~$\omega$.

\item
\textit{First order}: $[[\D,a],b] = 0$ for $a,b \in \A$.

\item
\textit{Orientability}:
There is a Hochschild $p$-cycle
\begin{subequations}
\label{eq:orient-cycle}
\begin{equation}
\cc = \sum_{\al=1}^n a_\al^0\ox a_\al^1\oxyox a_\al^p \in Z_p(\A,\A)
\label{eq:orient-cycle-Hoch}
\end{equation}
which is represented on~$\H$ by~$\Ga$, that is,
\begin{equation}
\pi_\D(\cc)
\equiv \sum_{\al=1}^n a_\al^0 \,[\D,a_\al^1] \dots [\D,a_\al^p] = \Ga.
\label{eq:orient-cycle-Gamma}
\end{equation}
\end{subequations}
\item
\textit{Closedness}:
$\Trw\bigl(\Ga\,[\D,a_1]\dots[\D,a_p]\,\<\D>^{-p}\bigr) = 0$ for all
$\row{a}{1}{p} \in \A$.

\item
\textit{Spin$^c$ or Morita equivalence:}
\label{cn:pmorita}
The $C^*$-module completion of $\H_\infty$ is a Morita equivalence
bimodule between $A$ and the norm completion of~$\CDA$.

\item 
\textit{Spin or Reality}:
There is an antiunitary operator $J$ on $\H$ such that
$J a^* J^{-1} = a$ for all $a \in \A$; $J^2 = \pm 1$;
$J \D J^{-1} = \pm\D$; and $J \Ga J^{-1} = \pm\Ga$ if $p$ is even;
with the same signs (depending only on $p \bmod 8$) as occur when 
$\A = \Coo(M)$ with $M$ a $p$-dimensional spin manifold, $\D$ is a
Dirac operator and $J$ is the charge conjugation.

\item 
\textit{Connectivity}:
\label{cn:conn}
There is an orthogonal family of projectors $p_j \in \A$ such that
$\Id_\A = \sum_j p_j$ and
$$
(a \in \A \word{with} [\D,a] = 0)  \iff a = \tsum_j \la_j p_j
\word{for some} \{\la_j\} \subset \C.
$$
\end{enumerate}

For a detailed discussion of how each of these conditions contributes
to the reconstruction of the smooth manifold $M$, we refer
to~\cite{Crux}. For instance, the orientability condition plays a role
in ensuring that all operators appearing under Dixmier traces above
are ``measurable'', which means that it does not matter which Dixmier
trace is used. Another point to notice is that, by the Serre--Swan 
theorem, there is a complex vector bundle $S \to M$ (that turns out 
to be the spinor bundle) for which $\H_\infty \subset \Ga(M,S)$, the 
latter being the $A$-module of continuous sections of~$S$.

Naturally, the conclusion that $M$ is a smooth manifold is a rather
strong one, and one may wonder whether some ``almost-manifolds'' could
also be allowed. For instance, one can build a spectral triple on an
orbifold, satisfying many, but not all, of the above postulates. While
one normally uses Lie groupoids~\cite{MoerdijkM} to describe
orbifolds, which entails noncommutative algebras, one could construct
a spectral triple on a quotient space $M/G$ using the algebra
$\A = \Coo(M)^G$ of smooth invariant functions, which satisfies many
of the properties listed above. We next show that this attempt fails,
precisely on account of the orientability condition.

\section{A $G$-invariant spectral triple}
\label{sec:orbi}

Consider the following example of a commutative spectral triple: let
$(M,g)$ be a compact boundaryless smooth Riemannian manifold with
metric $g$ of dimension~$p$, with a spin structure (so that in
particular $M$ is orientable), and let $G$ be a finite group acting by
isometries on~$M$. We shall suppose that $M$ is connected; if not, it
has finitely many components and we may restrict our attention to the
subgroup of~$G$ that preserves a given component. Assume that this
action lifts to an action of $G$ on the spinor bundle $S \to M$, for
which the Dirac operator $\Dslash_{\!g}$ is $G$-invariant. Now let
$\A := \Coo(M)^G$ be the algebra of $G$-invariant smooth functions
on~$M$, $\H := L^2(M,S)^G$ be the Hilbert space of $G$-invariant
$L^2$-spinors, and let $\D$ be the restriction of $\Dslash_{\!g}$ to a
selfadjoint operator on~$\H$. Then $(\A,\H,\D)$ is a spectral triple
over a commutative unital algebra.

(It might happen that an isometric action of $G$ on~$M$ lifts only 
to that of an extension $\widetilde G$ of~$G$ by~$\Z_2$, acting by 
automorphisms of the spinor bundle. Since $\widetilde G$ is also a 
finite group, the arguments below are not materially affected by this 
variant; so we may as well assume that the isometric action of $G$ 
lifts directly to the spinor bundle.)

We claim that this example \textit{cannot satisfy the orientation
condition} unless $G$ acts freely on~$M$.

If the action of $G$ on $M$ is free, then the orbit space $X = M/G$ is
a smooth Riemannian manifold and this spectral triple may be
identified with a Dirac spectral triple over $\Coo(X)$. In particular,
all the defining conditions listed in~\cite{Crux} will hold
automatically for $(\A,\H,\D)$. However, if the action is not free,
then $M/G$ is instead an orbifold. We do not need to consider this
quotient space as such, but much of the standard terminology of
orbifold theory is useful. We follow the notation of~\cite{MoerdijkM},
for the most part.

The isotropy subgroup of $x \in M$ is
$G_x := \set{h \in G : h\.x = x} \leq G$, and the fixed set of 
$h \in G$ is $\Sigma_h := \set{x \in M : h\.x = x} \subset M$. The
\textit{singular locus} is
$$
\Sigma_G := \bigcup_{h\neq 1} \Sigma_h = \set{x \in M : G_x \neq 1},
$$
so that $G$ acts freely on its complement. The exponential map at
$x \in M$ coming from the metric~$g$ gives a diffeomorphism 
$\exp_x : U \to W$ from a small open neighbourhood of $0 \in T_x M$ to
a neighbourhood $W$ of $x$ in~$M$. If $h \in G$, then $dh_x$ is an 
orthogonal linear transformation from $T_xM$ to $T_{h\.x} M$, such that
\begin{subequations}
\label{eq:exp}
\begin{equation} 
\exp_{h\.x} \circ dh_x = h \circ \exp_x : U \to h\.W,
\label{eq:exp-x-hx}
\end{equation}
When $h \in G_x$, we may take $W$ to be $G_x$-invariant and may assume
that each $dh_x$ preserves~$U$, so that
\begin{equation} 
\exp_x \circ dh_x = h \circ \exp_x : U \to W
\label{eq:exp-at-x}
\end{equation}
\end{subequations}
in that case. Note \cite[Lemma 2.10]{MoerdijkM} that if $h \in G_x$,
then $\set{y \in \Sigma_h : dh_y = 1 \text{ on } T_yM}$ is both closed
and open in $M$ by~\eqref{eq:exp-x-hx}, so that if $M$ is connected,
then $h \mapsto dh_x : G_x \to \Aut(T_x M)$ is injective for all
$x \in M$. Together with \eqref{eq:exp-at-x}, this shows that
$\Sigma_h$ has empty interior if $h \neq 1$. Thus the singular locus
$\Sigma_G$ is a closed subset of~$M$ with empty interior.

By transposition, each $dh_x$ acts on the cotangent spaces, taking
$T^*_{h\.x} M$ to $T^*_x M$, again orthogonally with respect to the 
transposed metric $g^{-1}$. In other words, $dh$ has a right action 
on covectors. 

Any $G$-invariant function $f$ equals its average over any cyclic 
subgroup of~$G$:
$$
f(x) = \frac{1}{|h|} \sum_{j=0}^{|h|-1} f(h^j\. x),
$$
where $|h|$ denotes the order of~$h$. If $h \in G_x$, taking
commutators with the Dirac operator gives
$$
[\D,f](x) = c(df(x))
= \frac{1}{|h|}\sum_{j=0}^{|h|-1} c(df(x) \rt dh_x^j),
$$
where $c$ denotes Clifford multiplication. Note that
\eqref{eq:exp-at-x} implies that the notation $dh_x^j$ is unambiguous:
$\exp_x \circ (dh_x)^j = h^j \circ \exp_x = \exp_x \circ (dh^j)_x$,
so that $(dh_x)^j = (dh^j)_x$ on~$T_x M$.

Let $x \in \Sigma_G$. The Cartan--Dieudonn\'e theorem, showing that
every orthogonal transformation on $T_x M$ is the product of at most
$p$ reflections, allows us to choose, for each $h \in G_x$, an
orthonormal basis for $T_x M$ such that
$$
dh_x =
\diag(R_{\theta_1},R_{\theta_2},\dots,R_{\theta_k},\pm 1,\dots,\pm 1),
$$
where each $R_{\theta_j}$ is a two-by-two rotation matrix (itself a
product of two reflections).

\vspace*{6pt}

Consider first the case where $G$ acts on $M$ by 
orientation-preserving isometries, so that each $dh_x$ is a rotation
in $SO(T_x M)$. Then we can write
$$
[\D,f](x) = c(df(x)) = c(u) + c(v_1) +\cdots+ c(v_k),
$$
where $u \in T^*_x M$ is invariant under $dh_x$ and
$v_r \rt dh_x = R_{-\theta_r} v_r$, for $r = 1,\dots,k$. Then
\begin{align*}
[\D,f](x) &= \frac{1}{|h|} \sum_{j=0}^{|h|-1} c(df(x) \rt dh_x^j)
\\
&= \frac{1}{|h|} \sum_{j=0}^{|h|-1}
c\bigl( u + R_{-\theta_1}^j v_1 +\cdots+ R_{-\theta_k}^j v_k \bigr)
\\
&= c(u) + \frac{1}{|h|} \sum_{r=1}^k \biggl( \sum_{j=0}^{|h|-1} 
R_{-\theta_r}^j \biggr) v_r = c(u),
\end{align*}
since $\sum_{j=0}^{|h|-1} R_{-\theta_r}^j = 0$: that follows since the
order of each rotation $R_{-\theta_r}$ divides~$|h|$.

Hence, if $f$ is any $G$-invariant function, $[\D,f](x)$ is fixed by
the rotation $dh_x$ and so it is (Clifford multiplication by) a
covector perpendicular to all the planes of rotation of~$h$.

Therefore, if $G_x \neq 1$ and if $a_0,a_1,\dots,a_p$ are 
\textit{$G$-invariant} functions, then $a_0\,[\D,a_1]\dots[\D,a_n]$
has vanishing skewsymmetrisation at~$x$:
\begin{equation} 
\frac{1}{p!} \sum_{\sg\in S_p} (-1)^\sg \sum_\al a_0(x)
\,[\D,a_{\sg(1)}](x) \dots [\D,a_{\sg(p)}](x) = 0,
\label{eq:skewsymm-null}
\end{equation}
since, under the usual linear isomorphism 
$\Cliff(T^*_x M, g_x) \simeq \La^\8 T^*_x M$, the left hand side
corresponds to a $p$-covector in $\La^p V$, where
$V = \bigcap_{h\in G_x} \ker({} \rt dh_x)$ has dimension less
than~$p$.

In particular, on inserting the entries of the Hochschild $p$-cycle
\eqref{eq:orient-cycle-Hoch} in Eqn.~\eqref{eq:skewsymm-null}, we
obtain
\begin{equation}
\Ga'(x) = \frac{1}{p!} \sum_{\sg\in S_p} (-1)^\sg \sum_\al a^0_\al(x)
\,[\D, a_\al^{\sg(1)}](x) \dots [\D, a_\al^{\sg(p)}](x) = 0
\label{eq:skew-Gamma}
\end{equation}
for all $x \in \Sigma_G$. We have written $\Ga'$ for the 
skewsymmetrized version of the represented cycle 
\eqref{eq:orient-cycle-Gamma}. Here $\Ga'$ may be regarded as an
endomorphism of the $\A$-module $\H_\infty$ and thereby as a
continuous section in $\Ga(M, \End S)$.

\vspace*{6pt}

Next consider the case where some $G_x$ contains a reflection.

Let $x \in \Sigma_G$ and $h \in G_x$ be such that $dh_x$ is a
reflection in the hyperplane $v^\perp$, for some unit covector
$v \in T^*_x M$. Thus, if $f$ is $G$-invariant, then
\begin{align*}
c(df(x)) &= [\D,f](x)
= \half \bigl( c(df(x)) + c(df(x) \rt dh_x) \bigr)
\\
&= \half c(df(x)) + \half\bigl( c(df(x)) - 2g_x(v, df(x))\,c(v) \bigr)
\\
&= c(df(x)) - g_x(v, df(x))\,c(v).
\end{align*}
Thus $g_x(v, df(x)) = 0$, so $df(x)$ is a covector in the hyperplane
$v^\perp$. Again, any $p$ such covectors have null exterior product,
so that $\Ga'(x) = 0$ in this case also.

\vspace*{6pt}

Suppose now that there are finitely many $G$-invariant functions
$\set{a_\al^j : j = 0,\dots,p,\ \al = 1,\dots,n}$ such that $\cc$ as
given by~\eqref{eq:orient-cycle-Hoch} is a Hochschild $p$-cycle, and
that Equation~\eqref{eq:orient-cycle-Gamma} holds with $\Ga$ being
proportional to Clifford multiplication by the volume element, as
happens when $\D$ is the Dirac operator. By using local orthonormal
bases of $1$-forms, one can check that the skewsymmetrization of this
section of the Clifford algebra bundle is, up to a constant, Clifford
multiplication by the Riemannian volume form $\vol_g$ on $M$. As such,
the skewsymmetrization is nowhere vanishing and the norm of
$\Ga'(y) \in \End S_y$ has a positive lower bound. Since we have shown
that for $G$-invariant functions there is always a neighbourhood of
the singular locus $\Sigma_G$ on which $\|\Ga'(y)\|$ is smaller than
any such lower bound, the orientability condition cannot be satisfied.

\vspace*{6pt}

The argument extends to any compact orbifold, since these are locally
of the form $M/G$ with $G$ finite. Indeed, such an orbifold $X$ is a
topological space that may be covered by finitely many orbifold charts
\cite{MoerdijkM}, each of the form $U/G$ where $U$ is an open
connected subset of~$\R^p$ and $G$ is a finite group of
diffeomorphisms of~$U$; one may impose a Riemannian metric on~$U$ for
which $G$ acts by isometries. Overlapping charts may use different
finite groups, but all are linked by compatibility homomorphisms. The
singular locus of the orbifold is the union of finitely many sets
$\Sigma_G$ in each such chart, and the above argument shows that the
candidate for a volume form, if it comes from a $G$-invariant $p$-form
on each~$U$, must always vanish on the singular locus, and therefore
cannot satisfy \eqref{eq:orient-cycle-Gamma}.

While the argument is local, it is not clear that for a general
orbifold we have a good candidate for a (commutative) algebra of
functions: this is highly dependent on the pasting
conditions~\cite{MoerdijkM}. When dealing with quotients $M/\Ga$,
where $\Ga$ is an infinite discrete group acting locally freely on~$M$
(that is, all isotropy subgroups are finite), there is such a
candidate algebra but no simple method of suitably encoding the local
freedom. It is therefore more appropriate to employ the groupoid
algebra in these situations.

\vspace*{6pt}

\subsubsection*{Acknowledgements}
We thank Gianni Landi for first asking us about the admissibility of
orbifolds as commutative geometries. Later, others have raised the
same question. This work was supported by an ARC grant, DP0211367, by
the Statens Naturvidenskabelige Forskningsr{\aa}d, Denmark. Support
from the University Complutense de Madrid and the Vicerrector\'{\i}a
de Investigaci\'on of the Universidad de Costa Rica is acknowledged.
JCV is grateful to Ryszard Nest for kind hospitality at the University
of Copenhagen.

\end{document}